\def\hrho{{\hat{\rho}}}
\def\hsig{{\hat{\sigma}}}
\def\sw#1{{\sb{(#1)}}}
\def\sco#1{{\sp{(\bar #1)}}} 
\def\su#1{{\sp{(#1)}}} 
\def\sul#1{{\sp{<#1>}}}
\def\br{{\bar\rho}}
\def\bs{{\bar\sigma}}

\def\proof{{\sl Proof.}\ }
\def\endproof{\hbox{$\sqcup$}\llap{\hbox{$\sqcap$}}}
\def\tens{\mathop{\otimes}}

\def\<{{\langle}}
\def\>{{\rangle}}

\def\id{{\rm id}} 
\def\eps{\epsilon}
\def\span{{\rm span}}
\def\q2{{q^{-2}}}

\def\cross{\mbox{$\times\!\rule{0.3pt}{1.1ex}\,$}}
\def\note#1{{}}
\def\eqn#1#2{\begin{equation}#2\label{#1}\end{equation}}

\def\Z{{\bf Z}}
\def\R{{\bf R}}
\documentstyle[12pt,epsf]{article}
\newtheorem{prop}{Proposition}[section]

\newtheorem{lemma}[prop]{Lemma}

\newtheorem{df}[prop]{Definition}
\newtheorem{ex}[prop]{Example}
\newtheorem{rem}[prop]{Remark}

\begin{document}

{\ }\qquad\qquad \hskip 3.4in DAMTP/96-28

{\ }\qquad\qquad \hskip 3.4in q-alg/9603016
~\vspace{2in}

\begin{center} {\large\bf CROSSED PRODUCTS BY A COALGEBRA}
\\ 
{\ }\vspace{25pt}\\
Tomasz Brzezi\'nski\footnote{Research
supported by the EPSRC grant GR/K02244}\vspace{21pt} \\
 Department of Applied Mathematics \& Theoretical Physics\\
University of Cambridge, Cambridge CB3 9EW\vspace{25pt}\\
Revised version, March 1997.
\end{center}
\vspace{63pt}
\begin{quote}\baselineskip 13pt
\noindent{\bf Abstract} We introduce the notion of a crossed product 
of an algebra by a coalgebra $C$, which  generalises the notion of a crossed
product by a bialgebra well-studied in the theory of Hopf algebras. The
result of such a crossed product 
is an algebra which is also a right $C$-comodule. We find the
necessary and sufficient conditions 
for two coalgebra crossed products be equivalent. 
We  show that the two-dimensional
quantum Euclidean group is a coalgebra crossed product. The paper is
completed with an appendix describing the  dualisation of construction
of coalgebra crossed products.

\bigskip

\end{quote}

\section{Introduction}
The notion of a crossed product by a bialgebra was first introduced
in the context of cohomology of algebras over a bialgebra
\cite{Swe:coh}. It was then generalised and throughly studied in the
relation to the Hopf-Galois theory of non-commutative rings and
invariant theory of algebraic groups \cite{DoiTak:cle}
\cite{Doi:equ} \cite{BlaCoh:cro} \cite{BlaMon:cro}. With the emergence
of quantum groups 
the theory of 
Hopf-Galois extensions and crossed products by Hopf algebras became
interesting also from geometric point of view (cf. \cite{Sch:hop} and
references therein). Hopf-Galois extensions 
are now understood as quantum group principal bundles
\cite{BrzMa:gau} \cite{Pfl:fib} 
and a certain kind of a crossed product, known as a cleft extension,
corresponds to a trivial quantum group principal bundle. 
 In search for
a suitable quantum group gauge theory on quantum homogeneous spaces 
we have recently proposed \cite{BrzMa:coa} a
generalisation of quantum group principal bundles in which the
structure quantum group is replaced by a coalgebra (a special case of
such a generalisation was already present in \cite{Sch:nor}, see in
particular Lemma~1.3). We also
introduced 
the notion of a trivial coalgebra principal bundle. As a vector space,
the total space of 
such a bundle is isomorphic to a tensor product of a base manifold
algebra and a structure coalgebra. It is therefore natural to expect
that as an algebra the total space is a certain kind of a
cleft extension of an algebra by a coalgebra. This leads naturally to
the notion of a crossed product by a coalgebra. In this paper we give
definition and study some properties of such crossed products.

Throughout the paper, all vector spaces are over a field $k$ of
generic characteristic. If not stated otherwise, by an algebra we mean
a unital associative algebra over $k$. The unit in an algebra is
denoted by 1 and the product by the juxtaposition of elements. We also
use $\mu$ to denote the product as a map. In a coalgebra $C$, $\Delta$
denotes the coproduct and 
$\eps :C\to k$ denotes the counit. We use the Sweedler notation to
denote the coproduct in $C$, $\Delta c = c\sw 1\tens c\sw 2$ (summation
understood), for any $c\in C$. By convolution product we mean a
product $*$ in a space of linear maps $C\to P$, where $P$ is an algebra,
given by $f*g(c) = f(c\sw 1)g(c\sw 2)$. A map $C\to P$ is said to be
convolution invertible if it is invertible with respect to this product.

\section{Crossed Products by a Coalgebra}
We begin by describing conditions which allow one to built an algebra
structure on a tensor product of an algebra and a vector space.
\begin{prop}
Let $M$ be an algebra, $V$ be a vector space and  $e\in V$. The vector
space $M\tens V$ is an algebra with unit $1\tens e$, and the product
such that
\begin{equation}
\forall x,y\in M \; \forall v\in V \qquad (x\tens e)(y\tens v) =
xy\tens v,
\label{leftlin}
\end{equation}
if and only if there exist linear maps $\hat{\sigma}:V\tens V\to
M\tens V$, $\hat{\rho}: V\tens M\to M\tens V$ which satisfy the
following conditions:\\
(a) $\forall x\in M\; \forall v\in V\quad \hat{\rho}(e,x) = x\tens e,
\quad \hrho(v,1) = 1\tens v$,\\
(b) $\hrho\circ(\id_V\tens\mu) =
(\mu\tens\id_V)\circ(\id_M\tens\hrho)\circ(\hrho\tens\id_M)$,\\
(c) $\forall v\in V\quad \hsig(e,v)=\hsig(v,e)=1\tens v$,\\
(d)
$(\mu\tens\id_V)\circ(\id_M\tens\hsig)\circ
(\hrho\tens\id_V)\circ(\id_V\tens\hsig) 
\! = \!(\mu\tens\id_V)\circ(\id_M\tens\hsig)\circ(\hsig\tens\id_V)$,\\
(e) $(\mu\tens\id_V)\circ(\id_M\tens\hsig)\circ
(\hrho\tens\id_V)\circ(\id_V\tens\hrho) \! 
= \! (\mu\tens\id_V)\circ(\id_M\tens\hrho)\circ(\hsig\tens\id_M)$,\\
where $\mu$ denotes product in $M$. The product $\mu_{M\tens V}$ in $M\tens V$
explicitly reads
\begin{equation}
\mu_{M\tens V} = (\mu^2\tens \id_V)\circ
(\id_M^2\tens\hsig)\circ(\id_M\tens\hrho\tens\id_V).
\label{prod.mv}
\end{equation}
\label{gen.cross.prod}
\end{prop}
\proof Let $\hsig: V\tens V\to M\tens V$, $\hrho: V\tens M\to M\tens
V$ be the linear maps that satisfy 
(a)-(e). We will show that the map $\mu_{M\tens V}$ given by
(\ref{prod.mv}) defines an algebra structure on $M\tens V$. To prove
that $\mu_{M\tens V}$ is associative we compute
\begin{eqnarray*}
\mu_{M\tens V}\circ(\id_M\tens\id_V\tens\mu_{M\tens
V})\!\!\!\!\!\!\!\!\!\!\!\!\!\!\!\!\!\!\!\!\!\!\!\!\!\!\!\!\!\!\!\!\!\!\!\!\!\!\!\!\!\!\!\!\!\!\!\!\!\!\!\!\!\!\!\!\!\!\!\!\!\!\!\!\!\!\!\!\!\!\!\!\!\!\!\!&&\\
 & = & (\mu^2\tens \id_V)\circ
(\id_M^2\tens\hsig)\circ(\id_M\tens\hrho\tens\id_V)\circ(\id_M\tens\id_V\tens\mu^2\tens
 \id_V)\circ\\
&&\circ
(\id_M\tens\id_V\tens\id_M^2\tens\hsig)\circ(\id_M\tens\id_V\tens\id_M\tens\hrho\tens\id_V)\\
&\stackrel{(b)}{=} & (\mu^2\tens \id_V)\circ
(\id_M^2\tens\hsig)\circ(\id_M\tens\mu^2\tens
 \id_V^2)\circ(\id_M^3\tens\hrho\tens\id_V)\circ \\
&&\!\!\!\!\circ(\id_M^2\tens\hrho\tens\id_M\tens\id_V)\circ(\id_M\tens\hrho\tens\id_M^2\tens\id_V)\circ
(\id_M\tens\id_V\tens\id_M^2\tens\hsig)\circ\\
&&\circ(\id_M\tens\id_V\tens\id_M\tens\hrho\tens\id_V)\\
& = & (\mu^4\tens \id_V) \circ (\id_M^4\tens\hsig) \circ
 (\id_M^3\tens\hrho\tens\id_V) \circ (\id_M^3\tens\id_V\tens\hsig) \circ\\
&&\circ (\id_M^2\tens\hrho\tens\id_V^2) \circ
 (\id_M^2\tens\id_V\tens\hrho\tens\id_V) \circ
 (\id_M\tens\hrho\tens\id_V\tens\id_M\tens\id_V)\\
&\stackrel{(d)}{=} & (\mu^4\tens \id_V) \circ (\id_M^4\tens\hsig) \circ
 (\id_M^3\tens\hsig\tens\id_V)\circ
(\id_M^2\tens\hrho\tens\id_V^2)\circ \\
&& \circ
 (\id_M^2\tens\id_V\tens\hrho\tens\id_V) \circ
 (\id_M\tens\hrho\tens\id_V\tens\id_M\tens\id_V)\\
&\stackrel{(e)}{=} & (\mu^4\tens \id_V) \circ (\id_M^4\tens\hsig) \circ
 (\id_M^3\tens\hrho\tens\id_V)\circ
(\id_M^2\tens\hsig\tens\id_M\tens\id_V)\circ \\
&& \circ
 (\id_M\tens\hrho\tens\id_V\tens\id_M\tens\id_V)\\
&=& (\mu^2\tens \id_V)\circ
(\id_M^2\tens\hsig)\circ(\id_M\tens\hrho\tens\id_V)\circ(\mu^2\tens
 \id_V\tens\id_M\tens\id_V)\circ\\
&&\circ
(\id_M^2\tens\hsig\tens\id_M\tens\id_V)\circ(\id_M\tens\hrho\tens\id_V\tens\id_M\tens\id_V)\\
&=&\mu_{M\tens V}\circ(\mu_{M\tens
V}\tens\id_M\tens\id_V). 
\end{eqnarray*}
To show that $1\tens e$ is a unit in $M\tens V$ we take any $x\in M$
and $v\in V$ and compute
\begin{eqnarray*}
\mu_{M\tens V}(1\tens e, x\tens v) & = & (\mu^2\tens \id_V)\circ
(\id_M^2\tens\hsig)\circ(1\tens\hrho(e,x)\tens v)\\
&\stackrel{(a)}{=} &
(\mu\tens\id_V)(x\tens\hsig(e,v))\stackrel{(c)}{=} x\tens v.
\end{eqnarray*}
Similarly one proves that $\mu_{M\tens V}(x\tens v, 1\tens e) = x\tens
v$. Therefore $\mu_{M\tens V}$ gives an associative algebra structure
on $M\tens V$ with unit $1\tens e$. To verify that property
(\ref{leftlin}) holds we take any $x,y\in M$ and $v\in V$ and compute
\begin{eqnarray*}
\mu_{M\tens V}(x\tens e, y\tens v) & = & (\mu^2\tens \id_V)\circ
(\id_M^2\tens\hsig)\circ(x\tens\hrho(e,y)\tens v)\\
&\stackrel{(a)}{=} &
(\mu^2\tens\id_V)(x\tens y\tens \hsig(e,v))\stackrel{(c)}{=} xy\tens v.
\end{eqnarray*}

Conversely assume that $M\tens V$ is an associative algebra with unit
$1\tens e$ and that (\ref{leftlin}) holds. Define maps $\hsig :V\tens
V\to M\tens V$ and $\hrho : V\tens M\to M\tens V$ by
$$
\hrho(v,x) = (1\tens v)(x\tens e), \qquad \hsig (v,w) = (1\tens
v)(1\tens w).
$$
Condition (\ref{leftlin}) implies that
\begin{equation}
(x\tens v)(y\tens e)= x\hrho(v,y), \qquad (x\tens
v)(1\tens w)= x\hsig (v,w),
\label{leftlin2}
\end{equation}
and relations (\ref{leftlin2}) imply that the product in $M\tens V$ has
the form (\ref{prod.mv}).

It remains to check that the maps $\hrho$ and $\hsig$ satisfy
conditions (a)-(e). Take any $x\in M$, $v\in V$. Then
$$
\hrho(v,1) = (1\tens v)(1\tens e) = 1\tens v ,\quad \hrho(e,x) =
(1\tens e)(x\tens e) = x\tens e,
$$
$$
\hsig(v,e) = (1\tens v)(1\tens e) = 1\tens v = (1\tens e)(1\tens v) =
\hsig (e,v).
$$
Therefore (a) and (c) hold. To check (b) take any $x,y\in M$ and
$v\in V$ and compute
\begin{eqnarray*}
\hrho(v,xy) &=& (1\tens v)(xy\tens e)
\stackrel{(\ref{leftlin})}{=}(1\tens v)(x\tens e)(y\tens e) =
\hrho(v,x)(y\tens e) \\
&\stackrel{(\ref{leftlin2})}{=}& (\mu\tens\id_V)\circ(\id_M\tens
\hrho)\circ(\hrho\tens \id_M)(v\tens x\tens y).
\end{eqnarray*}
To prove (d) take any $x\in M$, $v,w\in V$ and compute
\begin{eqnarray*}
(1\tens v)\!\!\!\!\!\!\!\!\!\!\!\!\!\!&&[(1\tens w)(x\tens e)] 
= (1\tens v)\hrho(w,x) \\
& =& (\mu^2\tens \id_V)\circ
(\id_M^2\tens\hsig)\circ(\id_M\tens\hrho\tens\id_V)(1\tens v\tens
\hrho(w,x)) \\
& = & (\mu\tens \id_V)\circ
(\id_M\tens\hsig)\circ(\hrho\tens\id_V)\circ(\id_V\tens
\hrho)(v\tens w\tens x).
\end{eqnarray*}
On the other hand the associativity of product in $M\tens V$ implies
that the above must be equal to
\begin{eqnarray*}
[(1\tens v)(1\tens w)](x\tens e) & = & \hsig(v,w)(x\tens e)\\
&\stackrel{(\ref{leftlin2})}{=}& 
(\mu\tens\id_V)\circ(\id_M\tens\hsig)\circ(\hsig\tens\id_V)(v\tens
w\tens x),
\end{eqnarray*}
Thus (e)
holds. Finally one verifies (d) considering
$$
(1\tens v)[(1\tens w)(1\tens u)] = [(1\tens v)(1\tens w)](1\tens u),
$$
for any $u,v,w\in V$. \endproof\medskip

This paper is concerned with a special case of the construction
described in Proposition~\ref{gen.cross.prod} in which $V=C$ is a
coalgebra equipped additionally with 
an entwining structure
\cite{BrzMa:coa}.
We say that  a coalgebra $C$  and an algebra $P$ are {\em
entwined}
if there is a map  $\psi :C\otimes P\to
P\otimes C$ such that
\eqn{ent.A}{\psi\circ(\id_C\tens
\mu)=(\mu\tens\id_C)\circ\psi_{23}\circ\psi_{12},\quad \psi(c\tens
1)=1\tens   c,\quad \forall c\in C}
\eqn{ent.B}{(\id_P\tens\Delta)\circ\psi=\psi_{12}\circ\psi_{23}
\circ(\Delta\tens\id_P),\quad (\id_P\tens\eps)\circ\psi=\eps\tens\id_P,}
where $\mu$ denotes multiplication in $P$, and
$\psi_{12} = \psi\tens\id_P$ and $\psi_{23}=\id_P\tens\psi$. 
We denote the action of $\psi$ on $c\otimes u\in C\otimes P$ by
$\psi(c\otimes u)= u_\alpha\tens c^\alpha$ (summation understood).

Furthermore we assume that there is a group-like $e\in C$,
i.e. $\Delta e=e\tens e$, $\eps(e) =1$ and a map $\psi^C:C\tens C\to
C\tens C$ such that for any $c\in C$
\begin{equation}
(\id\tens\Delta)\circ\psi\sp C=\psi\sp C_{12}\circ\psi^C_{23}
\circ(\Delta\tens\id),
\label{psiC.condition1}
\end{equation}
\begin{equation}
(\id\tens\eps)\circ\psi^C=\eps\tens\id,
\qquad \psi^C(e\otimes c) = \Delta c,
\label{psiC.condition2}
\end{equation}
where $\psi_{12}^C = \psi^C\tens\id_C$ and
$\psi_{23}=\id_C\tens\psi^C$. We denote the 
action of $\psi^C$ on $b\tens c$ by $\psi(b\tens c) = c_A\tens b^A$
(summation understood).

With these assumptions $P$ is a right $C$-comodule with a coaction
$\Delta_R u = \psi(e\tens u)$. Moreover the fixed point subspace
$M=P^{coC}_e =\{ x\in P | \Delta_R x = x\tens e\}$ is a subalgebra of
$P$. We call $(P, C,\psi, e,\psi^C)$ the {\em entwining data}. A
number of examples of entwining data may be found in
\cite{BrzMa:coa}. The 
main object of studies of this paper is contained in the following:
\begin{prop}
Let $(P, C,\psi, e,\psi^C)$ be entwining data and  let
$M=P^{coC}_e$. Assume that there are linear maps $\sigma :C\tens C\to
M$ and $\rho:C\tens P\to P$ such that for all $x,y\in M$,  $c\in C$:

(i) $\rho(e,x) = x, \quad \rho(c,1) = \eps(c)$;

(ii) $\rho (c\sw 1,x_\alpha)\tens {c\sw 2}^\alpha \in M\tens C$;

(iii) $\rho(c\sw 1,(xy)_\alpha)\tens c\sw 2^\alpha = \rho(c\sw
1,x_\alpha)\rho(c\sw 2^\alpha\sw 1,y_\beta)\tens c\sw 2^\alpha\sw
2^\beta$;

(iv) $\sigma(e, c) = \eps(c), \quad \sigma(c\sw 1, e_A)\tens {c\sw
2}^A = 1\tens c$.

Define a linear map $M\tens C \tens M\tens C\to M\tens C$, denoted by
juxtaposition, by
\begin{equation}
(x\tens b)(y\tens c) = x\rho(b\sw 1,
y_\alpha)\sigma(b\sw{2}^{\alpha}\sw 1, c_A)\tens
b\sw{2}^{\alpha}\sw 2^{A}.
\label{crossed.product.def}
\end{equation}
Then for any $x\in M$ and $c\in C$ $(1\tens e)(x\tens c) =
(x\tens c)(1\tens e) = x\tens c$.
The vector space $M\tens C$ is an associative algebra with product
(\ref{crossed.product.def}) if and only if
\begin{eqnarray}
\rho(a\sw 1, \sigma(b\sw 1, c_A)_\alpha)\sigma(a\sw 2^{\alpha}\sw
1, b\sw 2^{A}_{~B})\tens a\sw 2^{\alpha}\sw 2^{B}\!\!\!\!\!\!&&\nonumber \\
&& \!\!\!\!\!\!\!\!\!\!\!\!\!\!\!\!\!\!\!\!\!\!\!\!\!\!\!\!\!\!\!\!\!
\!\!\!\!\!\!  =
\sigma(a\sw 1, b_A)\sigma(a\sw 2^{A}\sw 1, c_B)\tens a\sw
2^{A}\sw 2^{B}
\label{cocycle}
\end{eqnarray}
and
\begin{eqnarray}
\rho(a\sw 1, \rho(b\sw 1,x_\alpha)_\beta)\sigma(a\sw 2^{\beta}\sw
1, {b\sw 2^{\alpha}}_{A})\tens a\sw 2^{\beta}\sw
2^{A}\!\!\!\!\!\!&&\nonumber \\ 
&& \!\!\!\!\!\!\!\!\!\!\!\!\!\!\!\!\!\!\!\!\!\!\!\!\!\!\!\!\!\!\!\!\!
\!\!\!\!\!\! =
\sigma(a\sw 1, b_A)\rho(a\sw 2^{A}\sw 1, x_\alpha)\tens a\sw
2^{A}\sw 2^{\alpha},
\label{twisted.module}
\end{eqnarray}
for any $a,b,c\in C$ and $x\in M$. This 
algebra is denoted by $M\cross_{\rho,\sigma}C$ and called
a {\em crossed product} by a coalgebra $C$. The pair $(\rho,\sigma)$ is called
the {\em crossed product data} for the entwining data $(P,C,\psi, e, \psi^C)$.
\label{crossed.product.prop}
\end{prop}
\proof This proposition can be viewed as a corollary of
Proposition~\ref{gen.cross.prod}. Namely given maps $\rho$ satisfying
 (ii) and  $\sigma$
 one can define maps $\hsig :C\tens
C\to M\tens C$ and $\hrho : C\tens M\to M\tens C$ by
$$
\hrho(x,c) = \rho(c\sw 1, x_\alpha)\tens c\sw 2^\alpha, \quad
 \hsig(b,c) = \sigma(b\sw 1,c_A)\tens b\sw 2^A.
$$
Conditions (i), (iii) and (iv) are precisely conditions (a)-(c) of
Proposition~\ref{gen.cross.prod}, and the map
(\ref{crossed.product.def}) is exactly as in (\ref{prod.mv}). Furthermore
conditions (\ref{cocycle}) and (\ref{twisted.module}) written in terms
of $\hrho$ and $\hsig$ are exactly the same as conditions (d) and (e)
of Proposition~\ref{gen.cross.prod}. By applying
Proposition~\ref{gen.cross.prod} to the present case the assertion
follows.  \endproof\medskip
\begin{rem}
\rm Note that while there is a one-to-one correspondence between maps
$\sigma$ and $\hsig$, different choices of $\rho$ and  may lead to the
same maps $\hrho$,  and thus to the same crossed
products. Therefore we will say that the crossed product data
$(\rho_1,\sigma)$ and $(\rho_2,\sigma)$ are {\em equivalent} if
$\hat{\rho_1} = \hat{\rho_2}$.
\label{rem.equivalent}
\end{rem}

Before we  specify to two special cases of crossed product
algebras we note that when $\psi(C\tens M)\subset M\tens C$ the
condition (iii) above 
is equivalent to
$$
\rho(c,xy) = \rho(c\sw 1,x_\alpha)\rho(c\sw 2^\alpha,y). 
$$
\begin{ex} \rm Let $C$ be a bialgebra $C=H$ and $P$ be a right
$H$-comodule algebra. Take $e=1$, and  the maps $\psi$ and
$\psi^C$  
$$
\psi(h\tens u) = u\sco 0\tens hu\sco 1, \qquad \psi^C(h\tens g) = g\sw
1\tens hg\sw 2,
$$
where $\Delta_R u = u\sco 0\tens u\sco 1$ (summation understood). $M$
is  the  fixed point subalgebra of  
a right $H$-comodule algebra $P$, and $\psi$ restricted to $M\tens C$
becomes the twist map, 
$\psi(h\tens x) = x\tens h$. At this point both $\psi$ and $P$ become
redundant in the setting of
Proposition~\ref{crossed.product.prop}. Conditions (i)-(iii) state that
$\rho$ is a weak normalised action of $H$ on $M$,
while
condition (iv) specifies the normalisation of $\sigma$, $\sigma(1,c)
= \sigma(c, 1)=\eps(c)$. Finally
 (\ref{cocycle}) states that $\sigma$ is a 2-cocycle, and
(\ref{twisted.module}) becomes equivalent to a twisted module
condition of \cite{DoiTak:cle}
\cite{BlaCoh:cro}. Therefore $M\cross_{\rho,\sigma} H$ is a
bialgebra crossed product.
\label{bialgebra.ex}
\end{ex}
\begin{ex} \rm Let
$C$ be a braided bialgebra $B$ and $P$ be 
a right braided $B$-comodule algebra both living in the same braided category,
and $e=1$ (for the review on braided bialgebras see
e.g. \cite{Ma:bey}). Then the maps $\psi$ and $\psi^C$ are 
$$
\psi(b\tens u) = \Psi(b\tens u\sco 0)u\sco 1, \qquad \psi^C(b\tens c)
= \Psi(b \tens c\sw 1)c\sw 2,
$$
where $\Psi$ denotes the braiding and $\Delta_R u = u\sco 0\tens u\sco
1$. In this case $\psi$ restricted to
$B\tens M$ coincides with $\Psi$.
Conditions (i)-(iii) state now that
$\rho$ is a normalised braided weak action of $B$ on $M$, while
conditions (iv) fix the normalisation of $\sigma$ as in
Example~\ref{bialgebra.ex}. Finally, 
using the  diagrammatic technique introduced in
\cite{Ma:bra}, in which  $\Psi=\epsfbox{braid.eps}$,  products are
denoted by
$\epsfbox{prodfrag.eps}$, and coproducts by
$\epsfbox{deltafrag.eps}$,  conditions (\ref{cocycle}) and 
(\ref{twisted.module}) come out respectively:
\[ \epsfbox{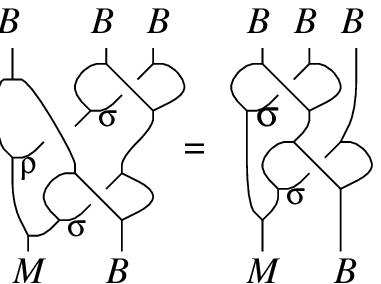}, \qquad ~ \qquad \epsfbox{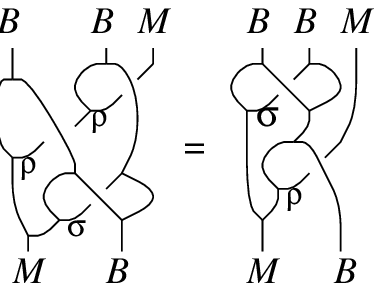}\]
Thus we obtain the generalisation of braided crossed products
introduced in \cite{Ma:cro} to the case of non-trivial $\sigma$. 
\label{braided.cross.product.ex}
\end{ex}

Therefore the notion of a crossed product by a
coalgebra in Proposition~\ref{crossed.product.prop} generalises both the notion
of a crossed product by a bialgebra and by a braided
bialgebra. The crossed product algebra $M\cross_{\rho,\sigma}C$ is a
left $M$-module with the natural action $x\tens y\tens c \mapsto
xy\tens c$. It is also a right
$C$-comodule with a coaction $\Delta_R : x\tens c \mapsto x\tens c\sw
1\tens c\sw 2$. Contrary to the  bialgebra crossed products
 $M\cross_{\rho,\sigma}C$ need not be a right $C$-comodule  algebra
even if $C$ is a  bialgebra. Instead, we have the following
\begin{lemma}
Let $M\cross_{\rho,\sigma}C$ be a crossed product as in
Proposition~\ref{crossed.product.prop} and let $\Delta_R
:M\cross_{\rho,\sigma}C \to M\cross_{\rho,\sigma}C \tens C$ be a natural
right coaction, $\Delta_R : x\tens c \to x\tens c\sw
1\tens c\sw 2$. Then
$$
\Delta_R((x\tens b)(y\tens c)) = (x\tens b\sw 1)(y_\alpha\tens
c_A)\tens b\sw 2^{\alpha A},
$$
for any $x,y\in M$ and $b,c\in C$, where the product on the right hand
side is given by (\ref{crossed.product.def}) even if $y_\alpha \not\in M$.
\label{twist.lemma}
\end{lemma}
\proof We compute
\begin{eqnarray*}
(x\tens b\sw 1)\!\!\!\!\!\!\!\!\!\!\!\!\!\!&&(y_\alpha\tens
c_A)\tens b\sw 2^{\alpha A} =   x\rho(b\sw
1,y_{\alpha\beta})\sigma(b\sw 2^\beta\sw 1, c_{AB})\tens b\sw
2^\beta\sw 2^B\tens b\sw 3^{\alpha A} \\
& \stackrel{(\ref{ent.B})}{=} &\!\!\!x\rho(b\sw
1,y_{\alpha})\sigma(b\sw 2^\alpha\sw 1, c_{A})\tens b\sw
2^\alpha\sw 2^A\sw 1\tens
b\sw 2^\alpha\sw 2^ A\sw 2 \\
& = &\!\!\! \Delta_R((x\tens b)(y\tens c)). \qquad \Box
\end{eqnarray*}

Now we give an example of a crossed product by a coalgebra coming
from the theory of coalgebra bundles \cite{BrzMa:coa}. 
\begin{ex} Let $P$, $C$, $\psi$, $e$, $\psi^C$ and $M$ be as in
Proposition~\ref{crossed.product.prop}. Furthermore assume that there
exists a convolution invertible map $\Phi : C\to P$ such that $\Phi(e)
=1$ and 
\begin{equation}
\psi\circ(\id_C\otimes\Phi) = (\Phi\otimes\id_C)\circ\psi^C.
\label{cov.phi}
\end{equation}
Define maps $\rho :C\tens P\to P$ and $\sigma :C\tens C\to P$ by
\begin{equation}
\rho(c,u) = \Phi(c\sw 1)u_\alpha\Phi^{-1}(c\sw 2^\alpha),
\qquad \sigma(b,c) = \Phi(b\sw 1)\Phi(c_A)\Phi^{-1}(b\sw 2^A).
\label{cleft.sigma}
\end{equation}
Then there is a crossed product
algebra $M\cross_{\rho, \sigma}C$. Explicitly the product in $M\cross_{\rho,
\sigma}C$ reads
\begin{equation}
(x\tens b)(y\tens c) = x\Phi(b\sw 1)y_\alpha \Phi(c_A)\Phi^{-1}(b\sw
2^{\alpha A}\sw 1)\tens b\sw 2^{\alpha A}\sw 2 .
\label{cleft.product}
\end{equation}
The algebra $M\cross_{\rho, \sigma}C$ is isomorphic to $P$. It is
called a {\em cleft extension} of $M$ by $C$ and is denoted 
by $M\cross_\Phi C$. 
\label{cleft.ex}
\end{ex}
\proof We first show that the output of the map $\sigma$
(\ref{cleft.sigma}) is indeed in $M$. We compute
\begin{eqnarray*}
\Delta_R\sigma(b,c) & \stackrel{(\ref{ent.A})}{=} & \Phi(b\sw
1)_\alpha \Phi(c_ A)_\beta\Phi^{-1}(b\sw 2^A)_\gamma \tens
e^{\alpha\beta\gamma}\\
 &\stackrel{(\ref{cov.phi})}{=} &\Phi(b\sw
1)\Phi(c_{AB})\Phi^{-1}(b\sw 3^A)_\gamma\tens b\sw 2^{B\gamma} \\
& \stackrel{(\ref{psiC.condition1})}{=} & \Phi(b\sw
1)\Phi(c_A)\Phi^{-1}(b\sw 2^A\sw 2)_\gamma\tens b\sw 2^A\sw 1^\gamma\\
& = &  \Phi(b\sw 1)\Phi(c_A)\Phi^{-1}(b\sw 2^A)\tens e .
\end{eqnarray*}
To derive the last equality we used the following property of
$\Phi^{-1}$,
\begin{equation}
\Phi^{-1}(c)\tens e = \Phi^{-1}(c\sw 2)_\alpha\tens c\sw 1^\alpha,
\label{cov.phi-1}
\end{equation}
which is easily obtained from (\ref{cov.phi}). Therefore $\sigma$ maps
$C\tens C$ to $M$ as required.

Next we  check that $\rho$ and $\sigma$ satisfy (i)-(iv) in
Proposition~\ref{crossed.product.prop}. The condition (i) is obvious. For
(ii) we take any $c\in C$ and $x\in M$ and compute
\begin{eqnarray*}
\Delta_R(\rho(c\sw 1 ,x_\alpha))\tens c\sw 2^\alpha & = & \psi(e, \Phi(c\sw
1)x_{\alpha\beta}\Phi^{-1}(c\sw 2^\beta)) \tens c\sw 3^\alpha \\
& \stackrel{(\ref{ent.B})}{=} &\psi(e, \Phi(c\sw
1)x_{\alpha}\Phi^{-1}(c\sw 2^\alpha\sw 1)) \tens c\sw 2^\alpha\sw 2 \\
& \stackrel{(\ref{ent.A})}{=} &\Phi(c\sw
1)_\beta x_{\alpha\gamma}\Phi^{-1}(c\sw 2^\alpha\sw 1)_\delta \tens
e^{\beta\gamma\delta}\tens c\sw 2^\alpha\sw 2 \\
& \stackrel{(\ref{cov.phi})}{=} &\Phi(c\sw
1) x_{\alpha\gamma}\Phi^{-1}(c\sw 3^\alpha\sw 1)_\delta \tens
c\sw 2^{\gamma\delta}\tens c\sw 3^\alpha\sw 2 \\
& \stackrel{(\ref{ent.B})}{=} &\Phi(c\sw
1) x_{\alpha}\Phi^{-1}(c\sw 2^\alpha\sw 2)_\delta \tens
c\sw 2^\alpha\sw 1^{\delta}\tens c\sw 2^\alpha\sw 3 \\
& \stackrel{(\ref{cov.phi-1})}{=} & \Phi(c\sw
1)x_{\alpha\beta}\Phi^{-1}(c\sw 2^\beta)\tens e\tens c\sw 3^\alpha .
\end{eqnarray*}

To check condition (iii) we take any $c\in C$, $x,y\in M$ and compute
\begin{eqnarray*}
\rho(c\sw 1 ,x_\alpha)\!\!\!\!\!\!\!\!\!\!\!\!\!\!\!&&
\rho(c\sw 2^\alpha\sw 1,y_\beta)\tens c\sw
2^\alpha\sw 2^\beta\\
 & = & \Phi(c\sw
1)x_{\alpha\beta}\Phi^{-1}(c\sw 2^\beta) \Phi(c\sw 3^\alpha\sw
1)y_{\gamma\delta}\Phi^{-1}(c\sw 3^\alpha\sw 2^\delta)\tens c\sw 3^\alpha\sw 3^\gamma \\
& \stackrel{(\ref{ent.B})}{=} & \Phi(c\sw
1)x_{\alpha}\Phi^{-1}(c\sw 2^\alpha\sw 1) \Phi(c\sw 2^\alpha\sw
2)y_\gamma\Phi^{-1}(c\sw 2^\alpha\sw 3^\gamma\sw 1) \tens c\sw
 2^\alpha\sw 3^\gamma\sw 2\\
& = &\Phi(c\sw
1)x_{\alpha}y_\gamma\Phi^{-1}(c\sw 2^{\alpha\gamma}\sw 1)\otimes c\sw
 2^{\alpha\gamma}\sw 2\\
&\stackrel{(\ref{ent.A})}{=}& \Phi(c\sw
1)(xy)_{\alpha}\Phi^{-1}(c\sw 2^\alpha\sw 1)\otimes c\sw 2^\alpha\sw 2\\
&\stackrel{(\ref{ent.B})}{=}&\Phi(c\sw
1)(xy)_{\alpha\beta}\Phi^{-1}(c\sw 2^\beta)\otimes c\sw 3^\alpha\\
&=& \rho(c\sw 1,(xy)_\alpha)\tens c\sw 2^\alpha.
\end{eqnarray*}
The first part of condition (iv) is obvious. For the second one we
have
\begin{eqnarray*}
\sigma(c\sw 1,e_A)\tens c\sw 2^A & = & \Phi(c\sw
1)\Phi(e_{AB})\Phi^{-1}(c\sw 2^B)\tens c\sw 3^A\\
&\stackrel{(\ref{psiC.condition1})}{=}& \Phi(c\sw
1)\Phi(e_{A})\Phi^{-1}(c\sw 2^A\sw 1)\tens c\sw 2^A\sw 2 \\
& \stackrel{(\ref{cov.phi})}{=} & \Phi(c\sw
1)\Phi(e)_\alpha\Phi^{-1}(c\sw 2^\alpha\sw 1)\tens c\sw 2^\alpha\sw 2\\
& = & \Phi(c\sw 1)\Phi^{-1}(c\sw 2)\tens c\sw 3\\
& = & 1\tens c.
\end{eqnarray*}
Therefore $\rho$ and $\sigma$ satisfy conditions (i)-(iv). To prove the
remaining part of the example we use the observation made in
\cite{BrzMa:coa} that $M\tens C$ is isomorphic to $P$ as a vector space
with the isomorphism $\Theta:x\tens c\to x\Phi(c)$ the inverse of
which is $\Theta^{-1}(u) = u_\alpha\Phi^{-1}(e^\alpha\sw 1)\tens
e^\alpha\sw 2$. It suffices to
show that $\Theta$ is an algebra isomorphism. We have
\begin{eqnarray*}
\Theta((x\tens b)(y\tens c)) & = & x\Phi(b\sw 1)y_\alpha \Phi(c_A)\Phi^{-1}(b\sw
2^{\alpha A}\sw 1)\Phi(b\sw 2^{\alpha A}\sw 2 )\\
& = & x\Phi(b\sw 1)y_\alpha \Phi(c_A)\eps(b\sw
2^{\alpha A})\\
& \stackrel{(\ref{ent.B},\ref{psiC.condition2})}{=}&
x\Phi(b)y\Phi(c) =  \Theta(x\tens b)\Theta(y\tens c).  
\end{eqnarray*}
Therefore $M\cross_{\Phi}C$ is isomorphic to $P$ as an algebra and
since $P$ is associative, so is $M\cross_{\Phi}C$. \endproof 

Since $P$ and $C$ are entwined by $\psi$, and $P$ is isomorphic to
$M\cross_{\Phi}C$ as an algebra, one would expect that the
cleft extension $M\cross_{\Phi}C$ and
$C$ are entwined by $\tilde{\psi} =
(\Theta^{-1}\tens\id_C)\circ\psi\circ(\id_C\tens\Theta)$. This is not
always the 
case, however. Instead we have
\begin{lemma}
$M\cross_{\Phi}C$ and $C$ are entwined by $\tilde{\psi}$ if and only if
$\psi^C (c\tens e) = e\tens c$, for any $c\in C$.
\end{lemma}
\proof Assume first that the hypothesis of the lemma is satisfied. Then
first of conditions (\ref{ent.A}) is satisfied by  $\tilde{\psi}$
since $\Theta$ is an algebra map. Explicitly, we have
\begin{eqnarray*}
\tilde{\psi}\circ(\id_C\tens\mu) & = &
(\Theta^{-1}\tens\id_C)\circ\psi\circ(\id_C\tens\Theta)\circ(\id_C\tens\mu) \\
& = &(\Theta^{-1}\tens\id_C)\circ\psi\circ
(\id_C\tens\mu)\circ(\id_C\tens\Theta\tens\Theta)\\ 
& = & (\Theta^{-1}\tens\id_C)\circ(\mu\tens\id_C)\circ
\psi_{23}\circ\psi_{12}\circ(\id_C\tens\Theta\tens\Theta)\\  
& = & (\mu\tens\id_C)\circ(\Theta^{-1}\tens\Theta^{-1}\tens\id_C)\circ
\psi_{23}\circ\psi_{12}\circ(\id_C\tens\Theta\tens\Theta)\\
& = & (\mu\tens\id_C)\circ
\tilde{\psi}_{23}\circ\tilde{\psi}_{12}.
\end{eqnarray*}
Here we used $\mu$ to denote the products in $P$ and in
$M\cross_{\Phi}C$. To prove the second of conditions (\ref{ent.A}) we first
notice, that 
$\tilde{\psi}$ can be explicitly written as
$$
\tilde{\psi}(b\tens x\tens c) =
x_{\alpha\beta}\Phi(c_{AB})\Phi^{-1}(e^{\beta B}\sw 1)\tens e^{\beta
B}\sw 2\tens b^{\alpha A}
$$
for any $b,c\in C$ and $x\in M$. In particular
\begin{eqnarray*}
\tilde{\psi}(b\tens 1\tens e) & = &
1_{\alpha\beta}\Phi(c_{AB})\Phi^{-1}(e^{\beta B}\sw 1)\tens e^{\beta
B}\sw 2\tens b^{\alpha A} \\
&=&\Phi(e_A\sw 1)\Phi^{-1}(e_A\sw 2)\tens e_A\sw 3 \tens b^{ A} =
1\tens e_A\tens b^A = 1\tens e\tens b,
\end{eqnarray*}
where we used the hypothesis to derive the last equality. The
verification of conditions (\ref{ent.B}) is easy.

Conversely, if we assume that $\tilde{\psi}$ entwines $M\cross_{\Phi}C$
with $C$ then the second of conditions (\ref{ent.A}) will imply
\begin{eqnarray*}
1\tens e\tens c & = &
1_{\alpha\beta}\Phi(c_{AB})\Phi^{-1}(e^{\beta B}\sw 1)\tens e^{\beta
B}\sw 2\tens c^{\alpha A} \\
&=&\Phi(e_A\sw 1)\Phi^{-1}(e_A\sw 2)\tens e_A\sw 3 \tens c^{ A}  = 1\tens
e_A\tens c^A,
\end{eqnarray*}
and thus the assertion follows. \endproof

\section{Equivalent Crossed Products}
\begin{prop}
Let  $M\cross_{\rho,\sigma}C$ be a  crossed product algebra as in 
Proposition~\ref{crossed.product.prop} associated to the entwining data
$(P,C,\psi, e,\psi^C)$, $M=P^{coC}_e$. Let $\gamma :C\to M$ be a convolution
invertible map such that $\gamma(e) =1$ and
\begin{equation}
\psi^C_{23}\circ\psi_{12}\circ
(\id_C\tens\gamma\tens\id_C)\circ(\id_C\tens\Delta) 
= (\gamma\tens\id_C\tens\id_C)\circ(\Delta\tens\id_C)\circ\psi^C.
\label{cov.gamma}
\end{equation}
Define the maps $\rho^\gamma:C\tens
P\to P$ and $\sigma^\gamma:C\tens C\to M$ by 
\begin{equation}
\rho^\gamma(c,u) =\gamma(c\sw 1)\rho(c\sw 2, u_\alpha)\gamma^{-1}(c\sw 3^\alpha),
\label{gauge.rho}
\end{equation}
\begin{equation}
\sigma^\gamma(b,c) = \gamma(b\sw 1)\rho(b\sw 2, \gamma(c_A\sw 1)_\alpha)\sigma
(b\sw 3^\alpha, c_A\sw 2)\gamma^{-1}(b\sw 4^A).
\label{gauge.sigma}
\end{equation}
Then $(\rho^\gamma,\sigma^\gamma)$ are the crossed product data and $M\cross_{\rho^\gamma,\sigma^\gamma}C$ is 
isomorphic 
to  $M\cross_{\rho,\sigma}C$ as an algebra. 
\label{coh.prop}
\end{prop}
\proof Clearly the output of $\sigma^\gamma$ is in $M$ by condition (ii) for $\rho$.
We now check that $\rho^\gamma$ and $\sigma^\gamma$ satisfy conditions (i)-(iv) of 
Proposition~\ref{crossed.product.prop}.
 Property (i) is immediate. Also, using (\ref{ent.B}) and the property
(ii) of $\rho$ we immediately obtain the property (ii) of $\rho^\gamma$. To verify
(iii) we take any $c\in C$, $x,y\in M$ and compute
\begin{eqnarray*}
\rho^\gamma(c\sw 1, (xy)_\alpha)\!\!\!\!\!\!\!\!\!\!\!\!\!\!\!&&\tens c\sw
2^\alpha =  \gamma(c\sw 1)\rho(c\sw 2,
(xy)_{\alpha\beta})\gamma^{-1}(c\sw 3^\beta)\tens c\sw 4^\alpha\\
&\stackrel{(\ref{ent.B})}{=}&\gamma(c\sw 1)\rho(c\sw 2,
(xy)_\alpha)\gamma^{-1}(c\sw 3^\alpha\sw 1)\tens c\sw 3^\alpha\sw 2\\
& = & \!\!\!\gamma(c\sw 1) \rho(c\sw 2,x_{\alpha})\rho(c\sw
3^\alpha\sw 1,y_\beta)\gamma^{-1}(c\sw 3^{\alpha}\sw 2^\beta\sw
1)\tens c\sw 3^{\alpha}\sw 2^\beta\sw 2\\
& = &\!\!\! \gamma(c\sw 1) \rho(c\sw 2,x_{\alpha})\gamma^{-1}(c\sw 3^\alpha\sw 1)\gamma(c\sw 3^\alpha\sw 2)\rho(c\sw
3^\alpha\sw 3,y_\beta)\\
&&\times\gamma^{-1}(c\sw 3^{\alpha}\sw 4^\beta\sw 1) 
\otimes c\sw 3^{\alpha}\sw 4^\beta\sw 2\\
& \stackrel{(\ref{ent.B})}{=} &\!\!\!\gamma(c\sw 1) \rho(c\sw 2,x_{\alpha\delta})\gamma^{-1}(c\sw 3^\delta)\gamma(c\sw 4^\alpha\sw 1)\rho(c\sw
4^\alpha\sw 2,y_\beta)\gamma^{-1}(c\sw 4^{\alpha}\sw 3^\beta\sw 1)\\ 
&& \otimes c\sw 4^{\alpha}\sw 3^\beta\sw 2\\
& = & \rho^\gamma(c\sw 1,x_\alpha)\gamma(c\sw 2^\alpha\sw 1)\rho(c\sw
2^\alpha\sw 2,y_\beta)\gamma^{-1}(c\sw 2^{\alpha}\sw 3^\beta\sw 1) 
\otimes c\sw 2^{\alpha}\sw 3^\beta\sw 2\\
& \stackrel{(\ref{ent.B})}{=} &\!\!\!\rho^\gamma(c\sw 1,x_\alpha)\gamma(c\sw 2^\alpha\sw 1)\rho(c\sw
2^\alpha\sw 2,y_{\beta\gamma})\gamma^{-1}(c\sw 2^{\alpha}\sw 3^\gamma) 
\otimes c\sw 2^{\alpha}\sw 4^\beta\\
& = &\!\!\!\rho^\gamma(c\sw 1,x_\alpha)\rho^\gamma(c\sw 2^\alpha\sw
1,y_\beta)\tens c\sw 2^\alpha\sw 2^\beta .
\end{eqnarray*}
In the third equality we used the fact that $\rho$ has the property
(iii).

To verify (iv) we take any $c\in C$ and compute
\begin{eqnarray*}
\sigma^\gamma(e,c) &=& \gamma(e)\rho(e,\gamma(c_A\sw
1)_\alpha)\sigma(e^\alpha,c_A\sw 2)\gamma^{-1}(e^A)\\
&=& \rho(e,\gamma(c_A\sw
1))\sigma(e,c_A\sw 2)\gamma^{-1}(e^A)= \gamma(c_A)\gamma^{-1}(e^A) =
\eps(c).
\end{eqnarray*}
Furthermore
\begin{eqnarray*}
\sigma^\gamma(c\sw 1, e_A)\tens c\sw 2^A\!\!\!\! & = &\!\!\!\! \gamma(c\sw 1)\rho(c\sw
2,\gamma(e_{AB}\sw 1)_\alpha)\sigma(c\sw 3^\alpha, e_{AB}\sw
2)\gamma^{-1}(c\sw 4^B)\tens  c\sw 5^A \\
& \stackrel{(\ref{ent.B})}{=} &\!\!\!\!\gamma(c\sw 1)\rho(c\sw
2,\gamma(e_{A}\sw 1)_\alpha)\sigma(c\sw 3^\alpha, e_{A}\sw
2)\gamma^{-1}(c\sw 4^A\sw 1)\tens\\
&&\tens c\sw 4^A\sw 2 \\
& \stackrel{(\ref{cov.gamma})}{=} &\!\!\!\!\gamma(c\sw 1)\rho(c\sw
2,\gamma(e)_{\beta\alpha})\sigma(c\sw 3^\alpha, e_A)\gamma^{-1}(c\sw
4^{\beta A}\sw 1)\tens c\sw 4^{\beta A}\sw 2 \\
& \stackrel{(\ref{ent.A})}{=} &\!\!\!\!\gamma(c\sw 1)\rho(c\sw
2,1)\sigma(c\sw 3, e_A)\gamma^{-1}(c\sw
4^{A}\sw 1)\tens c\sw 4^{A}\sw 2 \\
&=&\!\!\!\!\gamma(c\sw 1)\gamma^{-1}(c\sw 2)\tens c\sw 3 = 1\tens c.
\end{eqnarray*}
Therefore $\rho^\gamma$ and $\sigma^\gamma$ satisfy condition (i)-(iv) of
Proposition~\ref{crossed.product.prop}. Therefore we can define a map
 $\mu :M\tens C\tens M\tens C\to M\tens C$ by
(\ref{crossed.product.def}) with $\rho$ replaced by 
$\rho^\gamma$ and $\sigma$ replaced by $\sigma^\gamma$.
We consider a linear map
$\Theta : M\tens C\to M\cross_{\rho,\sigma}C$ given by 
\begin{equation}
\Theta(x\tens c) = x\gamma(c\sw 1)\tens c\sw 2.
\label{theta}
\end{equation}
It is clearly an isomorphism of vector spaces. We will show that for
all $x,y\in M$, $b,c\in C$, $\Theta(\mu(x\tens b\tens y\tens c)) =
\Theta(x\tens b)\Theta(y\tens c)$. To do so we first compute the map $\mu$
\begin{eqnarray*}
\mu(x\tens b\tens y\tens
c)
\!\!\!\!\!\!\!\!\!\!\!\!\!\!\!\!\!\!\!\!\!\!\!\!\!\!\!\!\!\!\!\!\!\!\!\!\!\!\!
\\
&=& x\gamma(b\sw 1)\rho(b\sw
2,y_{\alpha\beta})\gamma^{-1}(b\sw 3^\beta)\gamma(b\sw 4^\alpha\sw
1)\rho(b\sw 4^\alpha\sw 2, \gamma(c_{AB}\sw 1)_\delta)\times\\
&&\times \sigma(b\sw 4^\alpha\sw 3^\delta, c_{AB}\sw
2)\gamma^{-1}(b\sw 4^\alpha\sw 4^B)\tens b\sw 4^\alpha\sw 5^A\\
& \stackrel{(\ref{ent.B})}{=} & x\gamma(b\sw 1)\rho(b\sw
2,y_{\alpha})\rho(b\sw 3^\alpha\sw 1, \gamma(c_{A}\sw 1)_\delta)
\sigma(b\sw 3^\alpha\sw 2^\delta, c_{A}\sw
2)
\gamma^{-1}(b\sw 3^\alpha\sw 3^A\sw 1)\\
&&\tens b\sw 3^\alpha\sw
3^A\sw 2\\
& \stackrel{(\ref{cov.gamma})}{=} & x\gamma(b\sw 1)\rho(b\sw
2,y_{\alpha})\rho(b\sw 3^\alpha\sw 1, \gamma(c\sw 1)_{\beta\delta})
\sigma(b\sw 3^\alpha\sw 2^\delta, c\sw 2_{A})
\gamma^{-1}(b\sw 3^\alpha\sw 3^{\beta A}\sw 1)\\
&&\tens b\sw 3^\alpha\sw 3^{\beta A}\sw 2\\
& \stackrel{(\ref{ent.B})}{=} & x\gamma(b\sw 1)\rho(b\sw
2,y_{\alpha})\rho(b\sw 3^\alpha\sw 1, \gamma(c\sw 1)_{\beta})
\sigma(b\sw 3^\alpha\sw 2^\beta\sw 1, c\sw 2_{A})\\
&&\times \gamma^{-1}(b\sw 3^\alpha\sw 2^{\beta}\sw 2^ A\sw 1)
\tens b\sw 3^\alpha\sw 2^{\beta}\sw 2^ A\sw 2\\
& \stackrel{(iii)}{=} &x\gamma(b\sw 1)\rho(b\sw
2,(y\gamma(c\sw 1))_\alpha)
\sigma(b\sw 3^{\alpha}\sw 1, c\sw 2_{A})\gamma^{-1}(b\sw 3^{\alpha}\sw 2^
A\sw 1)\\
&&\tens b\sw 3^\alpha\sw 2 ^A\sw 2 .
\end{eqnarray*}
Thus
\begin{eqnarray*}
\Theta(\mu(x\tens b\tens y\tens c))\!\! & = &\!\!  x\gamma(b\sw 1)\rho(b\sw
2,(y\gamma(c\sw 1))_\alpha)
\sigma(b\sw 3^{\alpha}\sw 1, c\sw 2_{A})\tens b\sw 3^{\alpha}\sw 2^
A \\
 & = & \Theta(x\tens b)\Theta(y\tens c).
\end{eqnarray*}
Therefore $M\tens C$ is a unital associative algebra with product
$\mu$ and unit $1\tens e$. This algebra is a crossed product
$M\cross_{\rho^\gamma,\sigma^\gamma}C$ and it is isomorphic to
$M\cross_{\rho,\sigma}C$ with isomorphism $\Theta$. \endproof\medskip

\begin{df}
\rm We 
say that the crossed product data
$(\rho,\sigma)$ and $(\rho',\sigma')$ corresponding to the entwining
data $(P,C,\psi,e,\psi^C)$ are {\em gauge equivalent} if
there exists a convolution invertible map $\gamma: C\to M$ satisfying
conditions of Proposition~\ref{coh.prop} such that $(\rho',\sigma')$
are equivalent to $(\rho^\gamma, \sigma^\gamma)$ (in the sense of
Remark~\ref{rem.equivalent}). 
\label{gauge.def}
\end{df}
Convolution invertible maps $\gamma:C\to M$ satisfying
(\ref{cov.gamma}) were considered in \cite{BrzMa:coa} in the context
of trivial coalgebra $\psi$-principal bundles. They were termed gauge
transformations, since they induce a change of trivialisation
in a trivial $\psi$-principal bundle or cleft
extension of Example~\ref{cleft.ex}. Hence the terminology for the
crossed product data. It was also shown in
\cite{BrzMa:coa} that the set of all gauge transformations forms
a group with respect to convolution product.

{}From another point of view Proposition~\ref{coh.prop} generalises
the notion 
of equivalence of bialgebra crossed products of \cite{Doi:equ}
 to the case
of coalgebra crossed products. Precisely, we say that the crossed
products $M\cross_{\rho,\sigma}C$ and $M\cross_{\rho',\sigma'}C$,
corresponding to the same entwining data $(P,C,\psi,e,\psi^C)$ are {\it
equivalent} if there exists an algebra and a left $M$-module
isomorphism $\Theta : M\cross_{\rho',\sigma'}C \to
M\cross_{\rho,\sigma}C$, such that
\begin{equation}
\psi_{23}^C\circ\psi_{12}\circ(\id_C\tens\tilde{\Theta}) =
(\tilde{\Theta}\tens\id_C)\circ\psi^C,
\label{cov.theta}
\end{equation}
where $\tilde{\Theta}(c) = \Theta(1\tens c)$. Notice, that in
particular (\ref{cov.theta}) implies that $\Theta$ is a right $C$-comodule
isomorphism, i.e. $\Delta_R\circ\Theta =
(\Theta\tens\id_C)\circ\Delta_R$, where $\Delta_R$ is as in
Lemma~\ref{twist.lemma}.  
\begin{prop}
Crossed
products $M\cross_{\rho,\sigma}C$ and $M\cross_{\rho',\sigma'}C$,
corresponding to the same entwining data $(P,C,\psi,e,\psi^C)$ are 
equivalent if and only if the crossed product data $(\rho,\sigma)$
and $(\rho',\sigma')$ are gauge equivalent.
\label{equivalent.prop}
\end{prop}
\proof If $(\rho,\sigma)$
and $(\rho',\sigma')$ are gauge equivalent then the map $\Theta$
(\ref{theta}) is an algebra and a left $M$-module
isomorphism. Furthermore $\tilde{\Theta}(c) = \gamma(c\sw 1)\tens c\sw
2 = (\gamma\tens\id_C)\circ\Delta$. Using (\ref{cov.gamma}) we then
immediately obtain (\ref{cov.theta}).

Conversely, if $\Theta : M\cross_{\rho',\sigma'}C \to
M\cross_{\rho,\sigma}C$ is an algebra and a left $M$-module
isomorphism such that (\ref{cov.theta}) holds then we define the map
$\gamma :C\to M$ by $\gamma =
(\id_M\tens\eps)\circ\tilde{\Theta}$. Using (\ref{cov.theta}) and the
properties of the coaction $\Delta_R$ we obtain
\begin{eqnarray*}
LHS & = &
\psi^C_{23}\circ\psi_{12}\circ(\id_C\tens\gamma\tens\id_C)\circ
(\id_C\tens\Delta)\\ 
& = &\psi^C_{23}\circ\psi_{12}\circ(\id_C\tens\id_M\tens\eps\tens\id_C)
\circ(\id_C\tens\tilde{\Theta}\tens\id_C)\circ(\id_C\tens\Delta) 
\\
& = &
\psi^C_{23}\circ\psi_{12}\circ(\id_C\tens\id_M\tens\eps\tens\id_C)
\circ(\id_C\tens\Delta_R)\circ(\id_C\tens\tilde{\Theta})\\
& = &
\psi^C_{23}\circ\psi_{12}\circ(\id_C\tens\tilde{\Theta}) = 
(\tilde{\Theta}\tens\id_C)\circ\psi^C.
\end{eqnarray*}
On the other hand
\begin{eqnarray*}
RHS & = & (\gamma\tens\id_C\tens\id_C)\circ(\Delta\tens\id_C)\circ\psi^C \\
& = &
(\id_M\tens\eps\tens\id_C\tens\id_C)\circ(\tilde{\Theta}\tens\id_C\tens\id_C)
\circ(\Delta\tens\id_C)\circ\psi^C \\
& = &
(\id_M\tens\eps\tens\id_C\tens\id_C)\circ(\Delta_R\tens\id_C)\circ(\tilde{\Theta}\tens\id_C)\circ\psi^C\\
& = &  (\tilde{\Theta}\tens\id_C)\circ\psi^C.
\end{eqnarray*}
Therefore (\ref{cov.gamma}) holds. Clearly, $\gamma$ is convolution
invertible and $\gamma(e) = 1$. Moreover $\Theta(x\tens c) =
x\gamma(c\sw 1)\tens c\sw 2$. 
Since $\Theta$ is an algebra map  the following
expressions 
$$
\Theta(1\tens b)\Theta(x\tens c) =  \gamma(b\sw 1)\rho(b\sw 2,
(x\gamma(c\sw 1))_\alpha)\sigma(b\sw 3^{\alpha}\sw 1, c\sw
2_A)\tens b\sw 3^{\alpha}\sw 2^A ,
$$
$$
\Theta((1\tens b)(x\tens c)) = \rho'(b\sw 1, x_\alpha)\sigma'(b\sw
2^\alpha\sw 1, c_A)\gamma(b\sw
2^\alpha\sw 2^A\sw 1) \tens b\sw
2^\alpha\sw 2^A\sw 2
$$
are equal to each other. Applying $\id_M\tens \Delta$ and
$\id_M\tens\gamma^{-1}\tens \id_C$ we thus deduce
that
\begin{eqnarray}
\rho'(b\sw 1, x_\alpha)\sigma'(b\sw 2^\alpha\sw 1
,c_A)\tens b\sw 2^\alpha\sw
2^A\!\!\!\!\!\!\!\!\!\!\!\!\!\!\!\!\!\!\!\!\!\!\!\!\!\!\!\!\!\!\!\!\!\!\!\!\!\!\!\!\!\!\!\!\!\!\!\!\!\!\!\!\!\!\!\!\!\!\!\!\!\!\!\!\!\!\!\!\!\!\!\!\!\!\!\!\!\!\!\!\!
&&\label{prime}\\ 
& = &\!\!\! \gamma(b\sw 1)\rho(b\sw 2,
(x\gamma(c\sw 1))_\alpha)\sigma(b\sw 3^{\alpha}\sw 1, c\sw
2_A)\gamma^{-1}(b\sw 3^{\alpha}\sw 2^A\sw 1 )\tens b\sw 3^{\alpha}\sw 2^A\sw 2.
\nonumber
\end{eqnarray}
Setting $c=e$ in (\ref{prime}) and using property (iv) of $\sigma$ and
$\sigma'$ we obtain that $\rho'$ is equivalent to $\rho^\gamma$, and
setting $x=1$ in 
(\ref{prime}) we obtain (\ref{gauge.sigma}).
Therefore $(\rho,\sigma)$
and $(\rho',\sigma')$  are gauge equivalent. \endproof\medskip

Because any $\sigma$ satisfying condition (\ref{cocycle})
may be viewed as a 
generalisation of a cocycle, the results of
Proposition~\ref{coh.prop}  and Proposition~\ref{equivalent.prop} lead
 to the non Abelian cohomology theory
generalising that of  \cite[Section 6.3]{Ma:book}. To define a
coboundary, however, one would need to assume the existence of a
trivial cocycle $\sigma_{\rm triv}(b,c) = \eps(b)\eps(c)$. This
imposes additional conditions on the entwining data.
\begin{lemma}
The sufficient and necessary conditions for the existence of the
trivial cocycle $\sigma_{\rm triv}$ in the entwining data
$(P,C,\psi,e,\psi^C)$, with any map $\rho$ satisfying (i)-(iii) in
Proposition~\ref{crossed.product.prop}  are
\begin{equation}
\eps(e_A)c^A = c, \qquad \eps(c_A)\eps({b^A}_B)a^B =
\eps(b_A)\eps(c_B)a^{AB},
\label{trivial.cocycle}
\end{equation}
for any $a,b,c\in C$
\label{trivial.cocycle.lemma}
\end{lemma}
\proof In an elementary way one finds that first of conditions
(\ref{trivial.cocycle}) ensures second of the properties (iv) of
$\sigma_{\rm triv}$, while the second part of (\ref{trivial.cocycle})
is needed for the cocycle property (\ref{cocycle}) of $\sigma_{\rm
triv}$. \endproof\medskip

If the conditions of Lemma~\ref{trivial.cocycle.lemma} are satisfied
one defines a coboundary as a cocycle gauge equivalent to the trivial one,
i.e.
$$
\sigma(b,c) = \gamma(b\sw 1)\rho(b\sw 2,\gamma(c_A))\gamma^{-1}(b\sw
3^A).
$$
For example, conditions of  Lemma~\ref{trivial.cocycle.lemma} are
satisfied for braided bialgebra crossed products of
Example~\ref{braided.cross.product.ex}, therefore the coboundary can
be defined and the corresponding cohomology theory developed in this
case.

\section{$E_q(2)$ as a Crossed Product}
In this section we give an explicit example of a crossed product algebra
as defined in Proposition~\ref{crossed.product.prop}, namely, we show that
the two-dimensional quantum Euclidean group $E_q(2)$ is a
cleft extension as in Example~\ref{cleft.ex}, $E_q(2) =
X_q\cross_{\Phi}C$, 
where $X_q$ is the quantum hyperboloid \cite{Sch:int} and $C$ is
spanned by group-like elements $c_p$, $p\in \Z$ (other examples of
crossed products can be found in \cite{Brz:euc}). 
The algebra $P=E_q(2)$ is
generated by the elements $v,v^{-1}, n,\overline{n}$ subject to the
relations \cite{VakKor:alg} \cite{Wor:unb}
$$
vn = q^2nv, \qquad v\overline{n} = q^2 \overline{n}v, \qquad n\overline{n}
=q^2 \overline{n}n \qquad vv^{-1} = v^{-1} v =1,
$$
where $q\in k^*$. One can define a Hopf algebra structure on $E_q(2)$,
but this is not 
important in our construction. Instead we 
define the entwining structure $\psi:C\tens P\to P\tens C$ by
$$
\psi(c_p\tens v) = v\tens c_{p+1}, 
\quad \psi(c_p\tens n) = n\tens
c_p + \mu q^{2p}v\tens c_p -\mu q^{2p} v\tens c_{p+1} ,
$$
$$
\psi(c_p\tens v^{-1}) = v^{-1}\tens c_{p-1}, \quad
 \psi(c_p\tens
\overline{n}) = \overline{n}\tens 
c_p + \nu q^{2p}v^{-1}\tens c_p -\nu q^{2p} v^{-1}\tens c_{p-1},
$$
where $\mu, \nu \in k^*$ or $\mu =\nu =0$. One extends $\psi$ to the
whole of $P$ 
requiring that (\ref{ent.A}) be satisfied. To prove that such an
extension is well-defined one needs to show that 
$$
\psi(c_p\tens (vn -q^2 nv)) =\psi(c_p\tens (v\overline{n} -q^2
\overline{n}v)) = 0
$$
$$
\psi(c_p\tens (n\overline{n} -q^2
\overline{n}n)) =  \psi(c_p\tens (vv^{-1} -1))=\psi(c_p\tens (v^{-1}v
-1)) = 0.
$$
This can be done by elementary computations. It is also clear that
(\ref{ent.B}) holds. Therefore $\psi$ defines the entwining structure
for $P$ and $C$. Next we choose $e=c_s$ and define the right coaction of
$C$ on $P$ by $\Delta_R(u) = \psi(c_s\tens u)$. When $\mu =\nu =0$, $P$
becomes a right $C$-comodule algebra provided $C$ is equipped with the
algebra structure of $k[c,c^{-1}]$ by $c^p = c_{p+s}$, and therefore we do not
discuss this case any further. Following \cite{BonCic:fre} one easily
finds that
the fixed point subalgebra $M$ of $P$ is generated by
$
z = v + \mu^{-1}q^{-2s} n$ and $\overline{z} = v^{-1}
+\nu^{-1}q^{-2s}\overline{n}.
$
The generators $z$, $\overline{z}$ satisfy the relation
$
z\overline{z} = q^2 \overline{z} z +(1-q^2),
$
and therefore $M$ is isomorphic to the quantum hyperboloid $X_q$.

To complete the entwining data we define $\psi^C :C\tens C\to C\tens
C$ by 
$$
\psi^C(c_p\tens c_r) = c_r\tens c_{p+r-s}.
$$ 
It clearly
satisfies conditions (\ref{psiC.condition1},\ref{psiC.condition2}).

\begin{prop}
Let $\Phi : C\to E_q(2)$ be a linear map given by
$\Phi(c_p) = v^{p-s}$. Then $\Phi$ is convolution invertible and it satisfies
(\ref{cov.phi}). Therefore $E_q(2)= X_q\cross_{\Phi} C$
with a trivial cocycle $\sigma (c_p,c_r) =1$ and the map $\rho
:C\tens E_q(2)\to E_q(2)$,  
$$
\rho(c_p, v^{\pm 1}) =1,
$$
$$
\rho(c_p, n) = q^{2p}( q^{-2s}n +
\mu (v-1)), \quad \rho(c_p, \overline{n}) = q^{2p}( q^{-2s}\overline{n} +
\nu (v^{-1}-1)).
$$
\end{prop}
\proof Clearly $\Phi$ is convolution invertible and $\Phi^{-1}(c_p) =
v^{-p+s}$. Also, $\Phi(e) = \Phi(c_s) = 1$. To check that $\Phi$
satisfies (\ref{cov.phi}) we compute 
\begin{eqnarray*}
\psi(c_p\tens \Phi(c_r)) & = & \psi(c_p\tens v^{r-s}) = v^{r-s} \tens c_{p+r-s}\\
& = &\Phi(c_r) \tens c_{p+r-s}
 = (\Phi\tens \id_C)\circ\psi^C(c_p\tens c_r).
\end{eqnarray*}
By Example~\ref{cleft.ex} $E_q(2)= X_q\cross_{\rho,\sigma} C$,
with $\rho$ and $\sigma $ as specified. \endproof
\medskip

Although one can  easily
define a Hopf algebra structure on $C$ (e.g. $C$
could be a group algebra of any group $G$ such that $\# G =
\#\Z$) one equally easily finds that $E_q(2)$ is never a right
$C$-comodule algebra except when $(\mu,\nu)=(0,0)$. Furthermore,
because $\psi(C\tens X_q)$ is {\it not} a subset of $X_q\tens C$ there
is no braiding making $E_q(2)$ a braided $C$-comodule
algebra. Therefore the notion of a coalgebra crossed product developed in this
paper is truly needed for description of $E_q(2)$.
\medskip

\begin{center}
\bf ACKNOWLEDGMENTS
\end{center}
 The studies of crossed products by a coalgebra
were triggered by the work on coalgebra gauge theory \cite{BrzMa:coa},
a research carried out in collaboration with Shahn Majid. I am
grateful to 
him for continuing discussions.

\appendix
\section*{Appendix. Dual Crossed Products}
\setcounter{section}{1}
\addtocounter{prop}{-1}
In \cite{BrzMa:coa} it was observed that the definition of the
entwining structure $P$, 
$C$, $\psi$ (eqs. (\ref{ent.A}-\ref{ent.B})) possess the following
self-duality property. If one interchanges $P$ with $C$, $\mu$ with
$\Delta$, the unit map $\lambda\mapsto\lambda.1$ with $\eps$ and
reverses the order of composition of maps then the set of conditions
(\ref{ent.A}-\ref{ent.B}) remains unchanged. One can use this
self-duality property to dualise the notion of a crossed product algebra and
thus obtain the dual crossed product coalgebra. In this section we
present the result of such a dualisation.

We begin with an algebra $P$ and a coalgebra $C$ entwined by
$\psi$. Assume that there is an algebra character $\kappa :P\to k$. As
shown in \cite{BrzMa:coa}, the map $(\kappa\tens\id)\circ\psi$ is a
right action of $P$ on $C$. Moreover the linear space $J\sb\kappa
=\span\{c^\alpha\kappa(u_\alpha)-c\kappa(u)|c\in C,\ 
u\in P\}$ is
a coideal. Hence $M=C/J_\kappa$ is a coalgebra. 
Finally we assume that there is a
linear map $\psi^P :P\tens P\to P\tens P$ such that for any $u\in P$
\begin{equation}
\psi^P\circ(\id\tens\mu) = (\mu\tens\id)\circ\psi^P_{23}\circ\psi_{12}^P,
\label{psiP.condition1}
\end{equation}
\begin{equation}
\psi^P(u\tens 1) = 1\tens u, \qquad (\kappa\tens\id)\circ\psi^P =
\mu .
\label{psiP.condition2}
\end{equation}
Notice that (\ref{psiP.condition1})-(\ref{psiP.condition2}) come from
conditions 
(\ref{psiC.condition1})-(\ref{psiC.condition2}) by the dualisation
procedure explained 
above with the map $k\to C$, $\lambda\mapsto \lambda e$ replaced by
$\kappa$. We denote $\psi^P(u\tens v) = v_A\tens u^A$ (summation
understood). We say that $(C,P,\psi,\kappa,\psi^P)$ are {\em dual
entwining data}. Dualising Proposition~\ref{crossed.product.prop} we
thus obtain
\begin{prop}
Let $(C,P,\psi,\kappa,\psi^P)$ be dual entwining data, $M=C/J_\kappa$
with a canonical surjection $\pi :C\to M$, and let $\br :C\to P\tens
C$ and $\bs :M\to P\tens P$ be linear maps. We will denote the action
of these maps by $\br(c) =c\sul 1\tens c\sul 2$ and $\bs(m) = m\su
1\tens m\su 2$ (summation understood). Assume that for all $j\in J\sb
\kappa$, $m\in M$, $u\in P$ and $c\in C$ we have:

(i') $\kappa(c\sul 1)\pi(c\sul 2) = \pi(c), \quad c\sul
1\eps(c\sul 2) = \eps (c)$;

(ii') $j\sul 1u_\alpha \tens j\sul 2^\alpha \in P\tens J_\kappa$;

(iii') 
\begin{eqnarray*}
c\sul 1u_\alpha\tens \pi(c\sul
 2^\alpha\sw 1)\!\!\!\!\!\!\!\!\!\!&& \tens \pi(c\sul 2^\alpha\sw 2)\\
&& =
c\sw 1\sul 1(c\sw 2\sul 1u_\alpha)_\beta\tens\pi(c\sw 1\sul
2^\beta)\tens\pi(c\sw 2\sul 2^\alpha);
\end{eqnarray*}

(iv') $\kappa(m\su 1)m\su 2 = \eps(m), \quad m\su 1u_A\kappa(m\su 2^A)
= \eps(m)u.$

Then the vector space $M\tens P$ is a coalgebra with a coproduct
\begin{equation}
\Delta(m\tens u) = \pi(c\sw 1)\tens c\sw 2\sul 1\left(\pi(\sw 3)\su
1u_A\right)_\alpha\tens \pi(c\sw 2\sul 2^\alpha)\tens\pi(c\sw 3)\su
2^A,
\label{crossed.coproduct}
\end{equation}
where $c\in\pi^{-1}(m)$, and a counit $\eps\tens\kappa$ if and only if
for any $u\in P$ and $c \in C$, $m=\pi(c)$
\begin{eqnarray}
&&c\sw 1\sul 1(\pi(c\sw 2)\su 1u_A)_\alpha \tens \pi(c\sw 1\sul
2^\alpha)\su 1 {\pi(c\sw 2)\sul 2^A}_B \tens \pi(c\sw 1\sul
2^\alpha)\su 2^B \nonumber \\
&&\;\;\;\;\;\;\;\; = m\sw 1\su 1(m\sw 2\su 1 u_A)_B \tens m\sw 1\su
2^B\tens m\sw 2\su 2^A,
\label{cycle}
\end{eqnarray}
and 
\begin{eqnarray}
&& c\sw 1\sul 1(\pi(c\sw 2)\su 1u_A)_\alpha \tens c\sw 1\sul
2^\alpha\sul 1 {\pi(c\sw 2)\su 2^A}_\beta \tens \pi(c\sw 1\sul
2^\alpha\sul 2^\beta) \nonumber \\
&&\;\;\;\;\;\;\;\; = \pi(c\sw 1)\su 1(c\sw 2\sul 1 u_\alpha)_A \tens \pi(c\sw 1)\su
2^A\tens c\sw 2\sul 2^\alpha .
\label{twisted.comodule}
\end{eqnarray}
We denote the resulting coalgebra by $M\cross^{\br,\bs}P$ and call it a
{\em dual crossed product} coalgebra. We call $(\br,\bs)$ the {\em dual
crossed product data}
\label{crossed.coproduct.prop}
\end{prop}

Proposition~\ref{crossed.coproduct.prop} follows from the dual version of
Proposition~\ref{gen.cross.prod}. The latter can be proven by dualising the
proof of Proposition~\ref{gen.cross.prod} and is stated as follows
\begin{prop}
Let $C$ be an algebra, $V$ be a vector space and  $\kappa\in V^*$. The vector
space $C\tens V$ is a coalgebra with counit $\eps\tens \kappa$, and
the coproduct $\Delta_{C\tens V}$
such that
$$
\Delta\tens\id_V =
(\id_C\tens\kappa\tens\id_C\tens\id_V)\circ\Delta_{C\tens V}
$$
if and only if there exist linear maps $\hat{\sigma}:C\tens V\to
V\tens V$, $\hat{\rho}: C\tens V\to V\tens C$ which satisfy the
following conditions:\\
(a) $(\kappa\tens\id_C)\circ\hrho = \id_C\tens\kappa, \qquad
(\id_V\tens\eps)\circ\hrho = \eps\tens\id_V,$\\
(b) $(\id_V\tens\Delta)\circ\hrho =
(\hrho\tens\id_C)\circ(\id_C\tens\hrho)\circ (\Delta\tens\id_V)$,\\
(c) $(\id_V\tens\kappa)\circ\hsig = (\kappa\tens\id_V)\circ\hsig =
\eps\tens\id_V ,$\\
(d)
$(\id_V\tens\hsig)\circ(\hrho\tens\id_V)\circ(\id_C\tens\hsig)\circ(\Delta\tens\id_V)\! = \!(\hsig\tens\id_V)\circ(\id_C\tens\hsig)\circ(\Delta\tens\id_V)$,\\
(e) $(\id_V\tens\hrho)\circ
(\hrho\tens\id_V)\circ(\id_C\tens\hsig)\circ(\Delta\tens\id_V) 
= (\hsig\tens\id_C)\circ(\id_C\tens\hrho)\circ(\Delta\tens\id_V)$,\\
The coproduct $\Delta_{C\tens V}$ in $C\tens V$
explicitly reads
$$
\Delta_{C\tens V} = (\id_C\tens\hrho\tens\id_V)\circ
(\id_C^2\tens\hsig)\circ(\Delta^2\tens \id_V).
$$
\label{gen.cross.coprod}
\end{prop}
In the case of the dual crossed product the maps $\hrho : M\tens P\to P\tens M$ and $\hsig : M\tens P \to P\tens
P$ are given by
$$
\hrho(m,u) = c\sul 1u_\alpha \tens \pi(c\sul 2^\alpha), \qquad   
\hsig(m,u) =m\su 1u_A\tens m\su 2^A,
$$
for any $m\in M$, $u\in P$ and $c\in \pi^{-1}(m)$. 
Note  that the map $\hrho$ and consequently  $\Delta$ (\ref{crossed.coproduct})
are well-defined by condition (ii') and that the surjection $\pi$
dualises the natural inclusion
$M\hookrightarrow P$ consequently omitted in Sections~2 and 3. 

Similarly as for the crossed product data we say that
dual crossed product data $(\br_1,\bs)$ and $(\br_2, \bs)$ are
{\it equivalent} to each other if $\hrho_1 = \hrho_2$.

We now give an example of a dual crossed product coalgebra which is obtained
by the dualisation of Example~\ref{cleft.ex}.

\begin{ex}
Let $(C,P,\psi,\kappa,\psi^P)$ be dual entwining data, $M=C/J_\kappa$
with a canonical surjection $\pi :C\to M$, and
let $\Phi :C\to P$
be a convolution invertible map such that $\kappa\circ\Phi = \eps$ and
\begin{equation}
(\id_P\tens\Phi)\circ\psi = \psi^P\circ(\Phi\tens\id_P).
\label{cov.phi.dual}
\end{equation}
Define the maps $\br :C\to P\tens C$ and $\bs :M \to P\tens P$ by
\begin{equation}
\br(c) = \Phi(c\sw 1)\Phi^{-1}(c\sw 3)_\alpha\tens c\sw 2^\alpha, \qquad 
\bs(m) = \Phi(b\sw 1)\Phi^{-1}(b\sw 3)_A\tens\Phi(b\sw 2)^A ,
\end{equation}
where $b\in\pi^{-1}(m)$. Then there is a dual crossed product coalgebra
$M\cross^{\br,\bs}P$. The coproduct in $M\cross^{\br,\bs}P$ reads 
explicitly 
\begin{equation}
\Delta(m\tens u) = \pi(c\sw 1)\tens \Phi(c\sw 2)(\Phi^{-1}(c\sw  5)u)_{A
\alpha}\tens\pi(c\sw 3^\alpha)\tens\Phi(c\sw 4)^A,
\label{cocleft.sigma}
\end{equation}
where $c\in\pi^{-1}(m)$, and  $M\cross^{\br,\bs}P \cong C$ as coalgebras. 
\label{cleft.dual.ex}
\end{ex}
\proof This example is a dualisation of Example~\ref{cleft.ex} and
hence it can be proved following the same steps as in the proof of
Example~\ref{cleft.ex}. Therefore we do not repeat the full proof in
here. We only show that the map $\bs$ is well defined and we
write the explicit form of a coalgebra isomorphism  $M\cross^{\br,\bs}P
\cong C$. To prove the former we denote the function on the right hand
side of definition (\ref{cocleft.sigma}) of $\bs$ by $\tilde{\sigma}$
and compute
\begin{eqnarray*}
\tilde{\sigma}(\kappa(u_\alpha)c^\alpha) 
&\stackrel{(\ref{ent.B})}{=}&  \kappa(u_{\alpha\beta\gamma})\Phi(c\sw
1^\gamma)\Phi^{-1}(c\sw 3^\alpha)_A\tens\Phi(c\sw 2^\beta)^A \\
&\stackrel{(\ref{cov.phi.dual})}{=} & \kappa(u_{\alpha\beta B})\Phi(c\sw
1)^B\Phi^{-1}(c\sw 3^\alpha)_A\tens\Phi(c\sw 2^\beta)^A \\
&\stackrel{(\ref{psiP.condition2})}{=}& \Phi(c\sw
1)u_{\alpha\beta}\Phi^{-1}(c\sw 3^\alpha)_A\tens\Phi(c\sw 2^\beta)^A \\
&\stackrel{(\ref{cov.phi.dual})}{=}& \Phi(c\sw
1)u_{\alpha B}\Phi^{-1}(c\sw 3^\alpha)_A\tens\Phi(c\sw 2)^{BA}\\
&\stackrel{(\ref{psiP.condition1})}{=} &\Phi(c\sw
1)(u_{\alpha}\Phi^{-1}(c\sw 3^\alpha))_A\tens\Phi(c\sw 2)^{A}\\
& = &
\kappa(u)\Phi(c\sw 1)\Phi^{-1}(c\sw 3)\tens\Phi(c\sw 2) \\
& = & \tilde{\sigma}(\kappa(u)c).
\end{eqnarray*}
To derive the penultimate equality we used the following property of
$\Phi^{-1}$,
\begin{equation}
\Phi^{-1}(c)\kappa(u) = u_\alpha\Phi^{-1}(c^\alpha),
\label{cov.phi-1.dual}
\end{equation}
which can be easily derived from (\ref{cov.phi.dual}) or alternatively
obtained from (\ref{cov.phi-1}) by dualisation. Therefore
$J\sb\kappa\subset\ker\tilde{\sigma}$ and the map $\bs$ is
well-defined.

The coalgebra isomorphism $\Theta :C\to M\cross^{\br,\bs}P$ and its
inverse are given by
\begin{equation}
\Theta(c) = \pi(c\sw 1)\tens\Phi(c\sw 2), \qquad \Theta^{-1}(\pi(c)\tens u)
= \kappa((\Phi^{-1}(c\sw 2)u)_\alpha)c\sw 1^\alpha
\label{theta.dual}
\end{equation}
To see that the map $\Theta^{-1}$ is well-defined we denote the right
hand side of definition (\ref{theta.dual}) of $\Theta^{-1}$ by
$\theta$ and compute
\begin{eqnarray*}
\theta(\kappa(v_\gamma)c^\gamma\tens u) &\stackrel{(\ref{ent.B})}{=}&
\kappa\left(v_{\gamma\beta} \left(\Phi^{-1}\left(c\sw 2^\gamma\right)
u\right)_\alpha\right)c\sw 1^{\beta\alpha}\\ 
&\stackrel{(\ref{ent.A})}{=}& \kappa\left(\left(v_{\beta}
\Phi^{-1}\left(c\sw 2^\beta\right) 
u\right)_\alpha\right)c\sw 1^{\alpha}\\
& \stackrel{(\ref{cov.phi-1.dual})}{=}&\kappa(v)\kappa\left(\left(\Phi^{-1}(c\sw
2)u\right)_\alpha\right)c\sw 1^\alpha\\
& = & \theta(\kappa(v)c\tens u).
\end{eqnarray*}
Therefore $J_\kappa\subset\ker\theta$ and $\Theta^{-1}$ is
well-defined as stated. The proof that $\Theta$ is a comodule map and
that $\Theta^{-1}$ is its inverse can be obtained by careful
dualisation of the corresponding parts of the
proof of Example~\ref{cleft.ex} and hence we do not write it here.
\endproof\medskip

A dual crossed product $M\cross^{\br,\bs}P$ is a left $M$-comodule
with the coaction $\Delta_L(m\tens u) = m\sw 1\tens m\sw 2\tens
u$. Moreover $M\cross^{\br,\bs}P$ is a right $P$-module with the
action $m\tens u\tens v\mapsto m\tens uv$. 
We say that dual crossed products $M\cross^{\br,\bs}P$ and
$M\cross^{\br',\bs'}P$ 
corresponding to dual entwining data $(C,P,\psi,\kappa,\psi^P)$
are  {\em equivalent} if there exists a coalgebra and a left
$M$-comodule isomorphism $\Theta:M\cross^{\br,\bs}P\to
M\cross^{\br',\bs'}P$
such that
$$
\psi^P\circ(\tilde{\Theta}\tens\id_P) =
(\id_P\tens\tilde{\Theta})\circ\psi_{12}\circ\psi^P_{23} ,
$$
where $\tilde{\Theta} = (\eps\tens\id_P)\circ\Theta$. By dualising
Propositions~\ref{coh.prop} and \ref{equivalent.prop} we obtain the
following
\begin{prop}
Let $(C,P,\psi,\kappa,\psi^P)$ be dual entwining data, $M=C/J_\kappa$
with a canonical surjection $\pi :C\to M$, and
let $(\br,\bs)$, $(\br',\bs')$ be dual crossed product data. Then the
following statements are equivalent:

(1) Crossed products $M\cross^{\br,\bs}P$ and $M\cross^{\br',\bs'}P$
are equivalent.

(2) There exists a convolution invertible map $\gamma:M\to P$,
with the properties  
$\kappa\circ\gamma = \eps$ and 
$$
\psi^P\circ(\mu\tens\id_P)\circ(\gamma\tens\id_P\tens\id_P) =
(\id_P\tens\mu)\circ(\id_P\tens\gamma\tens\id_P)\circ\psi_{12}\circ\psi^P_{23},
$$
such that $(\br',\bs')$ are equivalent to
$$
\br^\gamma(c) = \gamma(\pi(c\sw 1))c\sw 2\sul 1\gamma^{-1}(\pi(c\sw
3))_\alpha\tens c\sw 2\sul 2^\alpha ,
$$
and
\begin{eqnarray*}
\bs^\gamma(\pi(c))\!\!\!\!\!\!\!\!\!\!\!\!\\
& = &\!\!\! \gamma(\pi(c\sw 1))c\sw 2\sul 1\pi(c\sw 3)\su
1_\alpha\gamma^{-1}(\pi(c\sw
4))_A\tens (\gamma(\pi(c\sw 2\sul 2^\alpha))\pi(c\sw 3)\su 2)^A,
\end{eqnarray*}
where $\br(c) = c\sul 1\tens c\sul 2$ and $\br(m) =m\su 1\tens m\su 2$.
\end{prop}

As before also this proposition can be proved by following the steps
of proofs of Propositions~\ref{coh.prop} and \ref{equivalent.prop},
therefore we skip the proof again. We only remark that the
map $\bs^\gamma$ is well defined by condition (ii').

In this way we were able to extend all the results of Sections 2 and
3 to the case of dual crossed products.


\end{document}